\documentclass[10pt,english]{amsart}
\usepackage[cp1251]{inputenc}
\usepackage[english]{babel}
\usepackage{amsmath}
\usepackage{amssymb}
\usepackage{amsfonts}
\usepackage[dvips]{graphicx}
\usepackage[colorlinks=true]{hyperref}
\newtheorem{theorem}{Theorem}

\newtheorem{defin}{Definition}

\newtheorem{corollary}{Corollary}

\newtheorem{example}{\sc Example}
\newtheorem{remark}{\sc Remark}

\newtheorem{nt}{\sc Note}

\newenvironment{pf}{\medskip\noindent{Proof:}
                \enspace}{\hfill \qed \newline \medskip}

\newenvironment{rmk}{\begin{remark}
      \em}{\end{remark}}

\begin{document}
\renewcommand{\refname}{References}
\thispagestyle{empty}
\title{On the Riemann-Hurwitz formula for   graph coverings}
\author{{A.D. Mednykh}}%
\address{Alexander Dmitrievich Mednykh
\newline\hphantom{iii} Sobolev Institute of Mathematics
\newline\hphantom{iii} Novosibirsk State University
\newline\hphantom{iii} 630090, Novosibirsk, Russia}%
\email{smedn@mail.ru}%

\thanks{\rm The work is supported by the Russian Foundation for Basic Research (grant 15–01–07906)}%

\maketitle {\small
\begin{quote}
\noindent{\sc Abstract. } The aim of  this paper is to present a few versions of the Riemann-Hurwitz formula for a regular branched covering of graphs. By a graph, we mean a finite connected  multigraph. The genus of a graph is defined as the rank of the first homology  group.  We consider a finite group acting on a graph, possibly with fixed and invertible edges, and the respective factor graph. Then, the obtained Riemann-Hurwitz formula relates genus of the graph with genus of the factor graph and   orders of the  vertex and edge stabilisers.

\vspace{8pt}

\vspace{4pt}
\end{quote}}

\section{Introduction}

Recall the classical Riemann-Hurwitz formula. Given surjective holomorphic map $\varphi: S \to  S^\prime $ between Riemann surfaces of genera $g $ and $g^\prime,$ respectively, one has
\begin{equation*}\label{RH_Classic}
2g-2={\rm deg}(\varphi)(2g^\prime-2)+\sum\limits_{ x\in S}(r_{\varphi}(x)-1),
\end{equation*}
where $r_{\varphi}(x)$ denotes the ramification index of $\varphi$ at $x.$ Let $G$ be a finite group of conformal automorphisms acting on $S$ and $\varphi: S \to  S^\prime=S/G $ is the canonical map induced by the group action. Then the above formula can be rewritten in the form
\begin{equation*}\label{RH_Regular}
2g-2=|G|(2g^\prime-2)+\sum\limits_{ x\in S}(|G^x|-1),
\end{equation*}
where $G^x$ stands for the stabiliser of  $x$ in $G$    and $|G^x|$ is the order of the stabiliser.
Remark that $S$ has only finite number of points with non-trivial stabiliser.

The latter formula has a natural discrete analogue.  Let $G$  be a finite group acting  on  the set of directed edges of a graph $X$ of genus $g$ freely and without invertible edges. Denote by $g^{\prime}$   genus of the factor graph $X^{\prime}=X/G.$ Then by \cite{BakerNorine} and \cite{Corry1} we have
\begin{equation*}\label{RH_Harmonic}
g-1=|G|(g^\prime-1)+\sum\limits_{ x\in V(X)}(|G^x|-1),
\end{equation*}
where $V(X)$   is the set of vertices of $X.$

The aim of this paper is to extend the above mentioned result to group actions with fixed and invertible edges. The main difficulty in this case is the correct definition of  a factor graph $X/G,$ when $G$ acts on a graph $X$ with invertible edges. There are at least three different ways to   define the graph $X/G.$ The first way is to consider  $X/G=(X/G)_{loop}$ as a graph with loops obtained as images of the invertible edges of $X$ under the canonical projection $X\to X/G.$ The second way is to consider the images of  invertible edges as semi-edges (or tails) of the factor graph $X/G=(X/G)_{tail},$ and the third one is to create the factor graph $X/G=(X/G)_{free}$ by deleting  loops (or semi-edges) that are the images of invertible edges. All the three ways are well known in the literature (\cite{MalNedSko}, \cite{Capo1}, \cite{BakerNorine}); they are effectively used in various questions of the graph theory.
\begin{figure}[h]
\center{\includegraphics[scale=0.93]{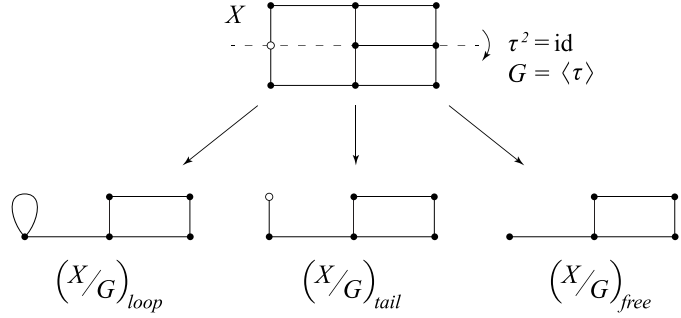}}
\caption{Different ways to define a factor graph $X/G.$}
\end{figure}%

\medskip
\section{Preliminary results and definitions}
In this section we introduce the notion of a graph with semi-edges that is a slightly more general then the standard  notion of a graph. This gives us a way to define the action of group on the graph with multiple edges and loops. Also, we are interesting in the group actions with fixed edges, as well as with invertible edges.  The factor space of such an action, in general, is not necessary a graph. But, it can be recognised as a graph with semi-edges.

Following \cite{MalNedSko} we define a  {\it graph\, with \,semi${}$-edges} as an ordered quadruple $X = (D,V;I,\lambda)$ where $D=D(X)$ is a set of ${\it darts,}$ $V=V(X)$ is a nonempty set of ${\it vertices,}$  which is required to be disjoint from $D,\, I$ is a mapping of $D $ onto $V,$ called the ${\it incidence\, function,}$ and $\lambda$ is an involutory permutation of $D,$ called the {\it dart-reversing} involution. For convenience or if $\lambda$ is not explicitly specified we sometimes write $\bar{x}$ instead of $\lambda x.$ Intuitively, the mapping $I$ assigns to each dart its ${\it initial \,vertex,}$ and the permutation $\lambda$ interchanges a dart and its reverse. The ${\it terminal\, vertex}$ of a dart $x$ is the initial vertex of $\lambda x.$  The $2$-orbits of $\lambda$ are called ${\it edges.}$ The $1$-orbits of $\lambda$ are called  {\it semi${}$-edges}
or {\it tails}. An edge is called   a ${\it loop}$   if $\lambda x \ne x$ and $I\lambda x = I x.$

We identify the set of edges $E(X)$  of $X$ with the following set of unordered pairs of darts:
$$E(X)=\{\{x, \bar{x}\}: x\in D(X),\,x\ne \bar{x}\} .$$
We will refer to the vertices $Ix$ and $I\bar{x}$ as {\it endpoints} of the edge $\{x, \bar{x}\}.$
In a similar way, the set of semi-edges $T(X)$  of $X$ is identified with the set
$$T(X)=\{\{x\}: x\in D(X),\,x=\bar{x}\} .$$

A {\it directed} edge of $X$ is an ordered pair $(x, \bar{x}),$ where $x\in D(X)\text{ and }x\ne \bar{x}.$ We note that all edges $\{x, \bar{x}\}\in E(X),$ including loops, are provided by  two directed edges $(x, \bar{x})$ and $(\bar{x}, x).$

A ${\it morphism \,of \,graphs} \,\,f:X=(D,V;I,\lambda)\to X^{\prime}=(D^{\prime}, V^{\prime};I^{\prime},\lambda^{\prime})$ is a function $f:D\cup V\to D^{\prime}\cup V^{\prime}$ such that $fD\subseteq D^{\prime},\,fV\subseteq V^{\prime},\,fI=I^{\prime}f \text{ and  }f\lambda = \lambda^{\prime}f.$ Thus, a morphism is an incidence-preserving mapping which takes vertices to vertices and edges to edges or semi-edges.  Note that the image of an edge can be an edge, a loop or a semi-edge, the image of a loop can be a loop or a semi-edge, and the image of a semi-edge can be just a semi-edge.

A bijective morphism $f:X\to X^{\prime}$ is called an {\it isomorphism,} and an isomorphism of $X$ onto itself  is called an {\it automorphism.} The group  ${\rm Aut}(X)$ of automorphisms of $X$  is a subgroup of $S_{D(X)}$  leaving invariant each of the sets $V(X),\, E(X),\,T(X)$ and preserving incidence.

We say that a group $G$ {\it acts} on $X$ if $G$ is a subgroup of ${\rm Aut}(X).$ We will refer to $X$ as a {\it graph} if it has the empty set of semi-edges
$T(X)=\emptyset.$

Let $X$ be a finite connected graph. We define the {\it genus} of $X$ to be the number
$$g(X)=1-|V(X)|+|E(X)|.$$ We note that $g(X)$ coincides with the Betti number of $X$ that is the rank of the first homology group ${\rm H}_1(X, \mathbb Z).$  Let $G$ be a finite group acting on the graph $X.$ An edge $e=\{x, \bar{x}\}\in E(X)$ is said to be {\it invertible} by $G$ if there is an element $g\in G$ such that $g$ sends $x$ to $\bar{x}$  and $\bar{x}$ to $x.$  An edge $e=\{x, \bar{x}\}\in E(X)$ is said to be {\it fixed} by $G$ if there is a non-trivial element $g\in G$  that fixes $x$ and $\bar{x}.$  We say that $G$ acts on $X$ {\it without invertible edges} if $X$ has no edges invertible by $G.$ Also,  $G$ acts on $X$ {\it without fixed edges} if $X$ has no edges fixed by $G.$

In this paper we will use two kinds of edge  stabiliser. The first one, $G^e,$ consists of
all elements of $G$ which fix the edge $e,$ that is fix both  $x$  and $\bar{x}.$ The second,  $G^{\{e\}},$ is  the setwise stabiliser of the set $e=\{x, \bar{x}\}$ in $G.$
We note that $|G^{\{e\}}| = 2|G^e|$ if the edge $e$ is invertible and $|G^{\{e\}}| = |G^e|$ otherwise.

\medskip
\section{Groups acting on a graph without invertible edges}
Our first result is the following theorem for groups acting on a graph without invertible edges.

\begin{theorem}\label{Main} Let $X$ be a graph of genus $g$ and $G$ is a finite group acting on $X$ without invertible edges. Denote by $g(X/G)$ genus of the factor graph  $X/G.$ Then

 \begin{equation*}\label{RH_WitoutEdgeRevising}
  g-1=|G|(g(X/G)-1)+\sum\limits_{ v\in V(X)}(|G^v|-1)-\sum\limits_{ e\in E(X)}(|G^e|-1),
 \end{equation*}
 where $V(X)$   is the set of vertices,   $E(X)$   is the set of edges of $X,\,\,G^x$ stands for the stabiliser of $ x\in V(X)\cup E(X)$ in $G$     and $|G^x|$ is the order of a stabiliser.
\end{theorem}
\begin{pf}Since $G$ acts on $X$ without invertible edges, the factor graph $X/G$ is well defined. The vertices of $X/G$ are orbits $Gv,\, v\in V(X),$ while the edges are orbits $Ge,\, e\in E(X).$ Vertices $Gv_1$ and $Gv_2$ are incident to an edge $Ge$ in $X/G$ if and only if $v_1$ and $v_2$ are incident to the edge $e$ in $X.$
Prescribe to every $\tilde{x}\in V(X/G)\cup E(X/G)$ a group $G_{\tilde{x}}$ isomorphic to $G^x,$ where $x$ is one of the preimages $\tilde{x}$ under the canonical map $\varphi:X\to X/G.$ Since $G$ acts transitively of fibres of $\varphi$ the group $G_{\tilde{x}}$ is well defined. One  can consider the graph  $X/G$ with prescribed groups  $G_v,\,v\in V(X/G)$ and $G_e,\,e\in E(X/G)$ as a graph of groups in sense of the Bass-Serre theory \cite{Bass}. We note that the fibre $\varphi^{- 1}(\tilde{x})$ of $\tilde{x}$ consists of $\frac{|G|}{|G_{\tilde{x}}|}$ elements. Hence,
\begin{equation}\label{VertexVolume}|V(X)|=\sum\limits_{v\in V(X)} 1=\sum\limits_{\tilde{v}\in V(X/G)}\frac{|G|}{|G_{\tilde{v}}|}\end{equation} and \begin{equation}\label{EdgeVolume}|E(X)|=\sum\limits_{e\in E(X)} 1=\sum\limits_{\tilde{e}\in E(X/G)}\frac{|G|}{|G_{\tilde{e}}|}.\end{equation}

By definition of genus from (\ref{VertexVolume}) and  (\ref{EdgeVolume}) we obtain

\begin{eqnarray}
g-1&=&|E(X)|-|V(X)|=  \sum\limits_{\bar{e}\in E(X/G)}\frac{|G|}{|G_{\tilde{e}}|}-\sum\limits_{\tilde{v}\in V(X/G)}\frac{|G|}{|G_{\tilde{v}}|} \nonumber\\
&=& |G|(\sum\limits_{\tilde{e}\in E(X/G)}1-\sum\limits_{\tilde{v}\in V(X/G)}1)\nonumber\\
&&+ \sum\limits_{\tilde{e}\in E(X/G)} \frac{|G|}{|G_{\tilde{e}}|}(1-|G_{\bar{e}}|)-\sum\limits_{\tilde{v}\in V(X/G)}\frac{|G|}{|G_{\tilde{v}}|}(1-|G_{\tilde{v}}|) \nonumber\\
&=& |G|(g(X/G)-1)+\sum\limits_{e\in E(X)}  (1-|G^{e}|)-\sum\limits_{v\in V(X)} (1-|G^{v}|)\nonumber\\
&=&|G|(g(X/G)-1)+ \sum\limits_{v\in V(X)} (|G^{v}|-1)-\sum\limits_{e\in E(X)}  (|G^{e}|-1).\nonumber
\end{eqnarray}
\end{pf}

 \section{Groups acting on a graph with invertible edges}
To describe the group action with edge revising we introduce a few definitions. First of all, we define the {\it barycentric subdivision} $X^{\prime}$ of a graph $X$ as a graph resulting from the subdivision of all edges in $X.$ The {\it subdivision} of some edge $e$ with endpoints $\{u,v\}$ yields a graph containing one new ({\it white}) vertex $w,$ and with an edge set replacing $e$ by two new edges, $\{u, w\}$ and $\{w, v\}.$ As a result, we consider the graph $X^{\prime}$ as a bipartite graph with black and white vertices that are in one to one correspondence with
vertices and edges of the graph $X,$ respectively.
The reverse operation, {\it smoothing} a vertex $w$ with regards to the pair of edges $(e,f)$ incident on $w,$ removes both edges containing $w$ and replaces $(e,f)$ with a new edge that connects the other endpoints of the pair. Here  we emphasise that only $2$-valent vertices can be smoothed.

Let now $G$ be a finite group acting on a graph $X,$ possibly with invertible edges. In this case, there are at least three different ways to define the factor graph $X/G.$

\vspace{8pt}
\noindent\textit{The factor graph with loops.} Define the image of an   edge $e$ with endpoints $\{u,v\}$ under the canonical map $X\to X/G$ to an edge $Ge$ with   endpoints $Gu $ and $Gv.$ If $e$ is an invertible edge  then the image of $e$ is a loop $Ge$ with the only one  endpoint $Gu = Gv.$ We denote   the obtained graph  $X/G$ by  $(X/G)_{loop}$ and its genus by $g(X/G)_{loop}.$
\smallskip

\noindent\textit{The factor graph with semi-edges.} Let $X^{\prime}$ be  the barycentric subdivision of the graph $X.$ For geometric  evidence, we identify the set $E(X^{\prime})$ of edges of $X^{\prime}$ with the set $D(X)$ of semi-edges of  $X.$ Then, the group $G$  naturally acts on the bipartite  graph $X^{\prime}$ preserving the vertex colour.  In particular, this means that $G$ acts on $X^{\prime}$  without invertible edges. That is the quotient $X^{\prime}/G$ is the well defined bipartite graph.
Note that the image of any invertible edge of $X$ in $X^{\prime}/G$ is now a bicolored edge with a  white vertex of valency one.

Consider  $(X/G)_{tail}$ as a graph obtained from the bipartite graph $X^{\prime}/G$ by smoothing  all $2$-valent white vertices. The images of   invertible edges are still bicolored edges with  white vertices of valency one. We will refer to them as {\it semi-edges} or the {\it tails} of graph  $(X/G)_{tail}.$ For the basic facts of   theory of graphs with semi-edges see Section 2 and the papers (\cite{MalNedSko}, \cite{Capo1}, \cite{Capo2}). Denote genus of $(X/G)_{tail}$ by $g(X/G)_{tail}.$

\vspace{8pt}
\noindent\textit{The factor graph without  semi-edges.} Denote by $(X/G)_{free}$ the graph obtained from $(X/G)_{tail}$ by removing all tails and replacing them by their black endpoint. Equivalently, $(X/G)_{free}$ can be obtained from $(X/G)_{loop}$ by removing all loops arising as the images of invertible edges. This kind of factor graphs was introduced by M. Baker and S. Norine in \cite{BakerNorine} for group $G$ generated by an involution. See also \cite{Corry2} for more detailed definitions. We write $g(X/G)_{free}$ for genus of $(X/G)_{free}.$

If the group $G$ acts on a graph $X$ with invertible edges then
$$g(X/G)_{free}=g(X/G)_{tail}\ne g(X/G)_{loop}.$$ If the action of $G$ has no invertible edges then all the three genera
coincide with $g(X/G).$

Our next result is  the following theorem.

\begin{theorem}\label{Semiedges}Let $X$ be a graph of genus $g$ and $G$ is a finite group acting on $X,$  possibly with invertible edges. Denote by $g(X/G)_{tail}$ genus of the factor graph  $(X/G)_{tail}.$ Then
$$g-1=|G|(g(X/G)_{tail}-1)+\sum\limits_{v\in V(X)}(|G^v|-1)-
\sum\limits_{e\in E(X)}(|G^e|-1)+\sum\limits_{ e\in E^{inv}(X)}|G^e|,$$
where $V(X)$   is the set of vertices,   $E(X)$   is the set of edges of $X,\,\,G^x$ is the stabiliser of $x\in V(X)\cup E(X)$  in $G,$   and $E^{inv}(X)$ is the set of invertible edges of $X.$
 \end{theorem}
\begin{pf}
Consider the set $D(X)$ of semi-edges of $X$. Let $e=(x,\bar{x})$, where $x\in D(X),$ be a directed edge of $X$.
Denote by $X^\prime$ the barycentric subdivision of graph $X$. Without loss of generality, one can assume that the edges of
$X^\prime$ are the elements of $D(X)$. Denote by $\partial_{1/2}e$ the vertex of $X^\prime$ subdividing the edge $e$ in two edges $x=h_e$ and $\bar{x}=h_{\bar{e}}$ of $X^\prime$.

 Since $G$ acts on $X^{\prime}$ without invertible edges, by Theorem~\ref{Main} we have
 \begin{equation}\label{RH_WithPrimes}
  g-1=|G|(g(X^{\prime}/G)-1)+\sum\limits_{ v\in V(X^{\prime})}(|G^v|-1)-\sum\limits_{ e\in E(X^{\prime})}(|G^e|-1).
 \end{equation}

 Here, the set of vertices $V(X^{\prime})$ of the bipartite graph $X^{\prime}$ is the union
 $V(X^{\prime})=B(X^{\prime})\cup W(X^{\prime})$ of the  sets of black and white vertices, where
 $$B(X^{\prime})=V(X)    \text{  and   } W(X^{\prime})=\{\partial_{1/2}e:\, e\in E(X) \}.  $$

 The stabiliser of a point $w=\partial_{1/2}e$ in $G$ consists of the elements of $G$ that permute the endpoints of $e$
 leaving $e$ invariant or the ones that fix $e.$ Hence, $G^{\partial_{1/2}e}=G^{\{e\}},$ where $G^{\{e\}}$ is the setwise stabiliser of the set $e=\{h_e,h_{\bar{e}}\}$  in $G.$ As a result, we obtain
\begin{equation}\label{RH_VerticesWithPrimes}
\begin{aligned}
  \sum\limits_{ v\in V(X^{\prime})}(|G^v|-1)&=\sum\limits_{ v\in V(X)}(|G^v|-1)+\sum\limits_{ e\in E(X)}(|G^{\partial_{1/2}e}|-1)\\
  &=\sum\limits_{ v\in V(X)}(|G^v|-1)+\sum\limits_{ e\in E(X)}(|G^{\{e\}}|-1).
\end{aligned}
\end{equation}

  For each $e\in E(X)$ we have $G^e=G^{h_e}=G^{h_{\bar{e}}}.$ Hence
\begin{equation}\label{RH_EdgesWithPrimes}
\begin{aligned}
  \sum\limits_{ e\in E(X^{\prime})}(|G^e|-1)&=\sum\limits_{ e\in E(X)}(|G^{h_e}|-1)+\sum\limits_{ e\in E(X)}(|G^{h_{\bar{e}}}|-1)\\
  &= 2\sum\limits_{ e\in E(X)}(|G^e|-1).
\end{aligned}
\end{equation}
  Subsituiting equations (\ref{RH_VerticesWithPrimes}) and (\ref{RH_EdgesWithPrimes}) into (\ref{RH_WithPrimes}) we obtain
  \begin{eqnarray}\label{RH_TwiceWithPrimes}
  g-1&=&|G|(g(X^{\prime}/G)-1)+\sum\limits_{ v\in V(X)}(|G^v|-1)-\sum\limits_{ e\in E(X)}(|G^e|-1)
  \nonumber\\&+&\sum\limits_{ e\in E(X)}(|G^{\{e\}}|-|G^e|).
 \end{eqnarray}
  Denote by $E^{inv}(X)$ the set of invertible edges of $X.$ Then $|G^{\{e\}}|=2|G^e|$ if $e\in E^{inv}(X)$ and
   $|G^{\{e\}}|=|G^e|$ otherwise. The smoothing of a white vertex in  graph $X^{\prime}/G$ decreases the number of vertices and the number of edges of the graph by one. So, it does not affect the genus $g(X^{\prime}/G)=1-|V(X^{\prime}/G)|+|E(X^{\prime}/G)|.$ Hence, by definition of $(X/G)_{tail}$ we have $g(X/G)_{tail}=g(X^{\prime}/G).$ Then (\ref{RH_TwiceWithPrimes}) can be rewritten in the form
$$
  g-1=|G|(g(X/G)_{tail}-1)+\sum\limits_{ v\in V(X)}(|G^v|-1)-\sum\limits_{ e\in E(X)}(|G^e|-1)
  +\sum\limits_{ e\in E^{inv}(X)}|G^e|.
$$
 \end{pf}%
\indent Since $g(X/G)_{free}=g(X/G)_{tail},$ as an immediate consequence of Theorem~\ref{Semiedges} we have the following statement.

\begin{theorem}\label{Free}Let $X$ be a graph of genus $g$ and $G$ is a finite group acting on $X,$  possibly with invertible edges. Denote by $g(X/G)_{free}$ genus of the factor graph  $(X/G)_{free}.$ Then

$$g-1=|G|(g(X/G)_{free}-1)+\sum\limits_{v\in V(X)}(|G^v|-1)-
\sum\limits_{e\in E(X)}(|G^e|-1)+\sum\limits_{ e\in E^{inv}(X)}|G^e|,$$
where $V(X)$   is the set of vertices,   $E(X)$   is the set of edges of $X,\,\,G^x$ is the stabiliser of $x\in V(X)\cup E(X)$  in $G,$   and $E^{inv}(X)$ is the set of invertible edges of $X.$
\end{theorem}

Following \cite{Corry2} we say that the group $G$ acts {\it harmonically} on a graph $X$ if $G$ acts freely on the set of darts $D(X)$ of $X$ or, equivalently, on the set of directed edges of $X$. In this case we have $|G^e|=1$ for each $e\in E(X).$ We have the following corollary from Theorem~\ref{Free}.

\begin{corollary}\label{Vertical}Let $X$ be a graph of genus $g$ and $G$ is a finite group acting on $X$ harmonically, possibly with invertible edges. Denote by $g(X/G)_{free}$ genus of the factor graph  $(X/G)_{free}.$ Then

$$g-1=|G|(g(X/G)_{free}-1)+\sum\limits_{v\in V(X)}(|G^v|-1)  +|E^{inv}(X)|,$$
where $V(X)$   is the set of vertices,   $E(X)$   is the set of edges of $X,\,\,G^v$ is the stabiliser of $v\in V(X)$  in $G,$   and $E^{inv}(X)$ is the set of invertible edges of $X.$
\end{corollary}

\begin{rmk}The Riemann-Hurwitz formula given in Corollary~\ref{Vertical}, up to notation, coincides with formula $(2.16)$ from \cite{BakerNorine}. Here, the set $E^{inv}(X)$ of invertible edges is exactly the set of {\it vertical} edges in terminology of \cite{BakerNorine}.
\end{rmk}

One more consequence of Theorem~\ref{Semiedges} is the following result.

\begin{theorem}\label{Loops}Let $X$ be a graph of genus $g$ and $G$ is a finite group acting on $X,$  possibly with invertible edges. Denote by $g(X/G)_{loop}$ genus of the factor graph  $(X/G)_{loop}.$ Then
 \begin{equation*}\label{RH_WithLoops}
  g-1=|G|(g(X/G)_{loop}-1)+\sum\limits_{ v\in V(X)}(|G^v|-1)-\sum\limits_{ e\in E(X)}(|G^{\{e\}}|-1),
 \end{equation*}
where $V(X)$ is the set of vertices, $E(X)$ is the set of edges of $X,\,\,G^v$ stands for the stabiliser of
$v\in V(X)$, $G^{\{e\}}$ stands for the stabiliser of the set consisting of two
semi-edges of $e\in E(X)$, and $|G^x|$ is the order of the stabiliser.
\end{theorem}
\begin{pf}
Let $X^\prime$ be the barycentric subdivision of $X$. Let $\varphi:X^\prime\to X^\prime/G$ be the canonical projection
and $\varphi:X\to (X/G)_{tail}$ be a map obtained from $\varphi^\prime$ by smoothing white $2$-valent vertices of
$X^\prime$ and $X^{\prime}/G$.
Denote by $T(X/G)$ the set of tails of $(X/G)_{tail}.$
The graph $(X/G)_{loop}$ can be obtained from $(X/G)_{tail}$ by replacing every tail of graph $(X/G)_{tail}$
with a loop.
Hence,
\begin{equation}\label{tails}
g(X/G)_{loop}=g(X/G)_{tail}+|T(X/G)|.
\end{equation}

Recall that each element of $T(X/G)$ is the image of an invertible edge $e\in E(X)$ under $\varphi.$
Hence, for any $\tilde{x}\in T(X/G)$ the fiber $\varphi^{-1}(\tilde{x})$ consists of $\frac{|G|}{|G^{\{e\}}|}$ invertible
edges of $X$. Since $G$ acts transitively on the fiber $\varphi^{-1}(\tilde{x})$, the number
$\frac{|G|}{|G^{\{e\}}|}$ does not depend on the choice of $e$ in the fiber. Therefore
\begin{equation}\label{tails2}
|T(X/G)|=\sum_{\tilde{x}\in T(X/G)}1=\sum_{e\in E^{inv}(X)}\frac{|G^{\{e\}}|}{|G|}
=\frac{1}{|G|}\sum_{e\in E^{inv}(X)}\big{|}G^{\{e\}}\big{|}.
\end{equation} By Theorem~\ref{Semiedges} we have
 \begin{equation}\label{eq_JJJ}
g-1=|G|(g(X/G)_{tail}-1)+\sum\limits_{v\in V(X)}(|G^v|-1)-
\sum\limits_{e\in E(X)}(|G^e|-1)+\sum\limits_{ e\in E^{inv}(X)}|G^e|.
\end{equation}
Since $|G^{\{e\}}|=2|G^e|$ if $e\in E^{inv}(X)$ and
   $|G^{\{e\}}|=|G^e|$ otherwise, from (\ref{tails}) and (\ref{tails2}) we get
\begin{equation}\label{eq_FFF}
|G|(g(X/G)_{tail}-1)=|G|(g(X/G)_{loop}-1)-2\sum_{e\in E^{inv}(X)}|G^{e}|.
\end{equation}
Also
\begin{equation}\label{eq_HHH}
\begin{aligned}
-\sum\limits_{e\in E(X)}(|G^e|-1)&=-\sum\limits_{e\in E(X)}(|G^{\{e\}}|-1)
+\sum\limits_{e\in E(X)}(|G^{\{e\}}|-|G^{e}|)\\
&=-\sum\limits_{e\in E(X)}(|G^{\{e\}}|-1)+\sum_{e\in E^{inv}(X)}|G^{e}|.
\end{aligned}
\end{equation}
Substituting (\ref{eq_FFF}) and (\ref{eq_HHH}) into (\ref{eq_JJJ}) we finally obtain
\begin{equation*}
  g-1=|G|(g(X/G)_{loop}-1)+\sum\limits_{ v\in V(X)}(|G^v|-1)-\sum\limits_{ e\in E(X)}(|G^{\{e\}}|-1).
 \end{equation*}
\end{pf}

\section{Acknowledgments} The author  is thankful to Gareth Jones, Tom Tucker, Toma\v{z} Pisanski and
Roman Nedela for fruitful discussions of the results of the paper.

\end{document}